\documentclass[preprint,review,authoryear]{elsarticle}

\usepackage{graphicx}
\usepackage{amssymb}
\usepackage{amsmath}
\usepackage{amsfonts}
\usepackage{graphicx}     

\journal{Control Engineering Practice}

\begin{document}

\begin{frontmatter}

\title{Pulse-Width Predictive Control for LTV Systems with Application to Spacecraft Rendezvous}

\author[aeroespacial-ESI]{R. Vazquez\corref{cor1}}

\author[aeroespacial-ESI]{F. Gavilan}
\ead{fgavilan@us.es}

\author[automatica-ESI]{E. F. Camacho}
\ead{eduardo@esi.us.es}

\address[aeroespacial-ESI]{Departamento de Ingenier\'{i}a Aeroespacial. Escuela Superior de Ingenieros. Universidad de Sevilla, Camino de los Descubrimientos s/n, 41092,Sevilla, Spain}
\address[automatica-ESI]{Departamento de Ingenier\'{i}a de Sistemas y Autom\'{a}tica. Escuela Superior de Ingenieros. Universidad de Sevilla, Camino de los Descubrimientos s/n, 41092,Sevilla, Spain}

\cortext[cor1]{Corresponding author, Fax number +34954486041, Email address: rvazquez1@us.es}

\begin{abstract}                
This work presents a model predictive controller (MPC) that is able to handle linear time-varying (LTV) plants with PWM control. The MPC is based on a planner that employs a PAM or impulsive approximation as a hot-start and then uses explicit linearization around successive PWM solutions for rapidly improving the solution by means of linear programming. As an example, the problem of rendezvous of spacecraft for eccentric target orbits is considered. The problem is modeled by the LTV Tschauner-Hempel equations, whose transition matrix is explicit; this is exploited by the algorithm for rapid convergence. The efficacy of the method is shown in a simulation study.\end{abstract}

\begin{keyword}
Spacecraft autonomy, Space robotics, Pulse-width modulation, Trajectory planning, Optimal trajectory, Linear Time-Varying Systems.
\end{keyword}

\end{frontmatter}

\section{Introduction}
Aerospace systems often need to be controlled by using pulse-width modulated (PWM) actuators, i.e., actuators whose output level is fixed and can only be turned on and off, such as spacecraft thrusters. It would be therefore desirable to use control design methods that directly take into account pulsed actuators. However, most feedback design and motion planning methods ignore variable width pulses and approximate the control variables either by impulses (which produce instantaneous changes in some combination of the states) or pulse-amplitude modulated (PAM) control. However, neither impulsive actuation nor PAM actuation capture with precision the behavior of pulsed actuators such as spacecraft thrusters. A more realistic model has to take into account that, typically, thrusters are ON-OFF actuators, i.e., the thrusters are not able to produce arbitrary forces, but instead can only be switched on (producing the maximum amount of force) or off (producing no force). These switching times are the only signals that can be controlled. This type of control signal is usually referred to as Pulse-Width Modulated (PWM). Control design with PWM actuation poses a challenge because the system becomes nonlinear in the switching times, even if the system is linear.

One can find in the literature several procedures to find  an equivalent PWM solution starting from a PAM solution (for instance in~\cite{Shieh1996673,Ieko1999123,Bernelli-Zazzera199864}). These methods allow to, given the PAM inputs of a system, compute  PWM inputs that produce a system output optimally approximating the output of the system when driven by the PAM signals. The results are based on the so-called Principle of Equivalent Areas, which computes the PWM signal so that it covers the same area as the PAM signal. However, while these procedures are quite effective in the sense that the output produced by the approximate PWM signals is very similar to the one produced by PAM signals, they assume that the plant is linear time-invariant.

In this paper, Model Predictive Control ({\sc MPC}) is used to directly find PWM signals to control the system. 
  {\sc MPC} (see, e.g.,~\cite{cabor:04})
is a family of methods that originated in  the late seventies  and has developed considerably
since then.   In {\sc MPC}, the process model is used to predict the future plant  
outputs, based on past and
current values and on the proposed optimal future control actions.  
These actions are calculated by the optimizer
taking into account the cost function as well as the
constraints. Since the plant is nonlinear in the control signals (ON-OFF times), the underlying optimization problem is nonlinear and possibly non-convex. To solve the problem, the algorithm starts from an initial guess computed by solving an optimal linear program with PAM or impulsive actuation, approximate the solution with ON-OFF thrusters, and then iteratively linearize around the obtained solutions to improve the PWM solution. While the idea of linearization to specifically compute optimal PWM control signals in the context of MPC is, to the best knowledge of the authors, original, it must be noted that local linearization techniques have been used for optimal trajectory problems in other contexts (see e.g.~\cite{sastry2002}).

As an application the problem of rendezvous of spacecraft is considered, i.e., the controlled close encounter of two space vehicles. Autonomous spacecraft rendezvous capabilities are becoming a necessity as access to space continues increasing. The field has become very active in recent years, with a rapidly growing literature. Among others, approaches based on trajectory planning and optimization~(\cite{How2008,arzelier2013,arzelier2011,louembet1,louembet2,louembet3,gaias1,gaias2}) and predictive control~(\cite{richards2003,rossi2002,asawa2006,Samara:2009,gavilan2012,larsson2006,Hartley2012695,lEOMANNI,scott1,scott2}) are emerging. 

Classically, in these approaches the problem of rendezvous is modeled by using impulsive maneuvers; one computes a sequence of (possibly optimal) impulses (usually referred to as $\Delta V$'s) to achieve rendezvous.

Recently,~\cite{IFAC_PWM, IFAC_PWM2}  introduced a  trajectory planning algorithm algorithm for spacecraft rendezvous that was able to incorporate PWM control signals. The former considered the linear time-invariant Clohessy-Wiltshire model (target orbiting in a \emph{circular} Keplerian  orbit, see~\cite{CW}). The latter extended the approach to \emph{elliptical} target orbits by using the  linear time-varying Tschauner-Hempel model (see~\cite{T-H}). Both methods start from an initial guess computed by solving an optimal linear program with PAM or impulsive actuation, approximate the solution with ON-OFF thrusters, and then iteratively linearize around the obtained solutions to improve the PWM solution. For both circular and elliptical target orbits the algorithms are simple and reasonably fast, and simulations favorably compare with an impulsive-only approach. These results were extended in~\cite{IFAC_PWM3} to a decreasing-horizon model predictive controller able to take into account orbital perturbations, disturbances or model errors.

In this paper, a \emph{receding} horizon model predictive controller with PWM inputs is formulated for general LTV plants and both alternatives (PAM or impulsive starting guess) are discussed in detail, with an application to rendezvous given at the end of the paper.

The structure of the paper is as follows. In Section~\ref{sect-model} the plant  model is introduced. Three types of inputs are considered: PWM, PAM and impulsive.
Section~\ref{sect-planning} follows with a formulation of the underlying optimization problem. Section~\ref{sect-method} describes a method that solves the planning problem using PWM signals. Section~\ref{sec-mpcrz} develops the model predictive controller. Next, Section~\ref{sec:rendezvous} describes the application to spacecraft rendezvous. Section~\ref{sect-simu} presents a simulation study of the method applied to spacecraft rendezvous. The paper finishes
with some remarks in Section~\ref{sect-rem}.

\section{System Model}
\label{sect-model} 
Consider a linear time-varying system given as
\begin{equation}
\dot x=A(t) x + B(t) u,\label{eqn-ltv}
\end{equation}
where $x\in \mathbb{R}^n$ is the state, $u\in \mathbb{R}^m$ is the input (control) vector, and $A(t)$ and $B(t)$ are, respectively, $n\times n$ and $n\times m$ matrices depending on time $t\geq 0$. 

Considering that, for some time $t_k\geq 0$, initial conditions $x(t_k)\in \mathbb{R}^n$ are given and the input is known, the solution to (\ref{eqn-ltv}) for $t>t_k$ is given by
\begin{equation}
x(t)=\Phi(t,t_k)x(t_k)+\int_{t_k}^t \Phi(t,s) B(s) u(s) ds,\label{eqn-sol}
\end{equation}
where $\Phi(t,t_k)$ is the \emph{system transition matrix}, see for instance~\cite{Rugh}. This matrix can be computed numerically (or analytically if possible) as the unique solution to the linear matrix differential equation 
\begin{eqnarray}
\dot \Phi(t,t_k)&=&A(t) \Phi(t,t_k),\quad t>t_k\\
\Phi(t_k,t_k)&=&\mathrm{I}.
\end{eqnarray}

To obtain an unified notation in term of the inputs, denote by $B_i(t)$ the $i$-th column of $B(t)$, corresponding to the $i$-th input $u_i(t)$, for $i=1,\hdots,m$. In the paper, time intervals starting at some initial time $t_k$ and ending at $t_{k+1}=t_k+T$ are considered, where $T$ will be an adequate sample time. Then equation (\ref{eqn-sol}) can be written as
\begin{eqnarray}
x(t_{k+1})&=&\Phi(t_{k+1},t_k)x(t_k)
+\sum_{i=1}^m \int_{t_k}^{t_{k+1}} \Phi(t_{k+1},s) B_i(s) u_i(s) ds,\label{eqn-sol2}
\end{eqnarray}

The objective is solving the problem with PWM inputs. In addition,  two other types of inputs are considered; they will be used as an intermediate step towards computing PWM inputs by the algorithm. All types of input are analyzed in the following sections.

\subsection{Pulse width-modulated (PWM) control}
\begin{figure}[t]
\centering
\includegraphics[width=0.5\columnwidth]{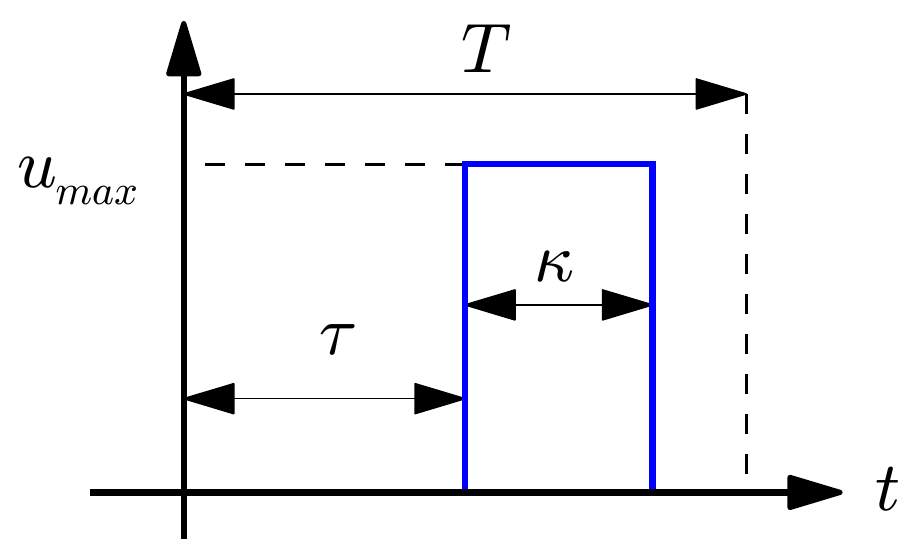}    
\caption{PWM Variables.}
\label{fig-pwm_variables}
\centering
\end{figure}

In the PWM case, each input $u_i$ is a pulse starting at time $\tau_{k,i}$ (relative to $t_k$) with pulse width  $\kappa_{k,i}$, with constant magnitude $u_{max}=u^W_{k,i}$, as shown in Fig.~\ref{fig-pwm_variables}, i.e.,
\begin{equation} \label{eq:control_PWM_def}
u_i\hspace{-1pt}(t)=\left\{ \begin{array}{ll} \vspace{1pt}
0, & t\in\left[t_k,t_k+\tau_{k,i}\right],\\\vspace{1pt}
u^W_{k,i}, & t\in\left[t_k+\tau_{k,i},t_k+\tau_{k,i}+\kappa_{k,i}\right],\\
0, & t\in\left[t_k+\tau_{k,i}+\kappa_{k,i},t_{k+1}\right],
 \end{array} \right.
\end{equation}
with $\kappa_{k,i}>0$, $\tau_{k,i}>0$ and $\tau_{k,i}+\kappa_{k,i}<T$, where the last constraint prevent the PWM signal to spill over to the next time interval.
Then, substituting $u_i(t)$ in (\ref{eqn-sol2}) one obtains
\begin{eqnarray}
x(t_{k+1})&=&\Phi(t_{k+1},t_k)x(t_k)
+\sum_{i=1}^m \left( \int_{t_k+\tau_{k,i}}^{t_k+\tau_{k,i}+\kappa_{k,i}} \hspace{-27pt}\Phi(t_{k+1},s) B_{i}(s) ds\right) u^W_{k,i} ,\,\quad\label{eqn-sol_pwm1}
\end{eqnarray}
and denoting 
\begin{equation}\label{eqn_BWdef}
B^W_{k,i}(\tau_{k,i},\kappa_{k,i})= \int_{t_k+\tau_{k,i}}^{t_k+\tau_{k,i}+\kappa_{k,i}} \Phi(t_{k+1},s) B_{i}(s) ds,\end{equation}
one can write the solution as
\begin{eqnarray}
x(t_{k+1})&=&\Phi(t_{k+1},t_k)x(t_k)
+\sum_{i=1}^m
B^W_{k,i}(\tau_{k,i},\kappa_{k,i})u^W_{k,i}.\qquad\label{eqn-sol_pwm}
\end{eqnarray}
There is an important difference between a PWM input and the PAM or impulsive inputs that will be subsequently introduced. While the latter can in principle take positive o negative values at different times, the former is fixed either as positive or negative for all time. Thus, typically a PWM model has twice number of inputs than a PAM/impulsive model. To make this explicit in the model (\ref{eqn-sol_pwm}), denote with a plus or minus super-index the positive or negative inputs, as follows
\begin{equation}
x(t_{k+1})
=\Phi(t_{k+1},t_k)x(t_k)
+\sum_{i=1}^m\left[
B^W_{k,i}(\tau_{k,i}^+,\kappa_{k,i}^+)u^{W+}_{k,i}-
B^W_{k,i}(\tau_{k,i}^-,\kappa_{k,i}^-)
u^{W-}_{k,i}\right],\qquad\label{eqn-sol_pwm}
\end{equation}
with $u^{W+}_{k,i},\tau_{k,i}^+,\kappa_{k,i}^+$ and $u^{W-}_{k,i},\tau_{k,i}^-,\kappa_{k,i}^-$ denoting, respectively, the magnitude, start, and width of the positive and negative $i$-th input pulses.

\subsection{Pulse amplitud-modulated (PAM) control}
In this case, each control $u_i(t)$ in (\ref{eqn-sol2}) is constant inside the interval $[t_k,t_{k+1}]$, and equal to $u^A_{k,i}$.
Then, substituting $u_i(t)$ in (\ref{eqn-sol2}) one obtains
\begin{eqnarray}
x(t_{k+1})&=&\Phi(t_{k+1},t_k)x(t_k)
+\sum_{i=1}^m \left( \int_{t_k}^{t_{k+1}} \hspace{-12pt}\Phi(t_{k+1},s) B_i(s) ds\right) u^A_{k,i} ,\quad\label{eqn-sol_imp}
\end{eqnarray}
and denoting by 
\begin{equation}\label{eqn-defBA}
B^A_{k,i}= \int_{t_k}^{t_{k+1}} \Phi(t_{k+1},s) B_i(s) ds,\end{equation}
 the solution  can be written as
\begin{eqnarray}
x(t_{k+1})&=&\Phi(t_{k+1},t_k)x(t_k)+\sum_{i=1}^m  B^A_{k,i} u^A_{k,i}.\label{eqn-sol_pam}
\end{eqnarray}

\subsection{Impulsive control}
 In this case, $u_i(t)=u^I_{k,i} \delta(t-(t_k+\tau_{k,i}))$, where $\delta(t)$ is Dirac's delta function, $t_k+\tau_{k,i}$ is the instant at which the impulse is given, and $u_{k,i}$ is the magnitude of the impulse. Then, assuming $0<\tau_{k,i}<T$ for all $i$ (all the impulses are given inside the considered time interval) and substituting $u_i(t)$ in (\ref{eqn-sol2}) one obtains
\begin{eqnarray}
x(t_{k+1})&=&\Phi(t_{k+1},t_k)x(t_k)
+\sum_{i=1}^m  \Phi(t_{k+1},t_k+\tau_i) B_i(t_k+\tau_i) u^I_{k,i} ,\label{eqn-sol_imp1}
\end{eqnarray}
and denoting by $B^I_{k,i}(\tau_{k,i})= \Phi(t_{k+1},\tau_i) B_i(t_k+\tau_i)$, 
\begin{eqnarray}
x(t_{k+1})&=&\Phi(t_{k+1},t_k)x(t_k)+\sum_{i=1}^mB^I_{k,i}(\tau_{k,i})u^I_{k,i}.\label{eqn-sol_imp}
\end{eqnarray}

\subsection{Discretization and compact notation}
Consider now a sequence of time instants $t_k=t_0+kT$, $k=0,\hdots$, and denote $x_k=x(t_k)$. Then, it is possible to write, for both PAM and impulsive control,
\begin{equation}
x_{k+1}=A_k x_k+B_k U_{k},\label{eqn-simp_evol}
\end{equation}
where $A_k=\Phi(t_{k+1},t_k)$, and $B_k$ and $U_{k}$ depend on the input type.
In the PWM case,  write 
\begin{equation}
x_{k+1}=A_k x_k+B^+_k U^+_k-B^-_k U^-_k,\label{eqn-simp_evol_PWM}
\end{equation}
where $B^+_k$ is a matrix whose $i$-th column is $B^W_{k,i}(\tau_{k,i}^+,\kappa_{k,i}^+)$ and $U^+_k$, $\tau_{k}^+$, and $\kappa_k^+$  are column vectors whose $i$-th entries are, respectively, $u^W_{k,i}$, $\tau_{k,i}^+$ and $\kappa^+_{k,i}$. The same definitions (with minus super-index) are used for the negative pulses. Then, to reach model (\ref{eqn-simp_evol}),  define
\begin{equation}
B_k = \left[ \begin{array}{cc}  B^+_k & -B^-_k
\end{array} \right],
\,
{U}_k = \left[ \begin{array}{c}  U^+_k \\
 U^-_{k}
\end{array} \right],
\,
{\tau}_k  = \left[ \begin{array}{c}\tau_{k}^+
\\Ê{\tau}_{k}^-
\end{array} \right],
\,
{ \kappa}_k = \left[ \begin{array}{c}  \kappa_k^+ \\
 \kappa_{k}^-
\end{array} \right].
\end{equation}

The definitions of $B_k$ are simpler for the other types of actuation. In the PAM case, $B_k$ is a matrix whose $i$-th column is $B^A_{k,i}$ and $U_k$ a column vector whose $i$-th entry is $u^A_{k,i}$. In the impulsive case, $B_k$ is a matrix whose $i$-th column is $B^I_{k,i}(\tau_{k,i})$, and $U_k$, $\tau_k$ are column vector whose $i$-th entries are, respectively, $u^I_{k,i}$ and $\tau_{k,i}$.

Next a compact formulation is developed to simplify the notation of the problem.
The state at time $t_{k+j+1}$, given the state $x_k$ at time $t_k$,  and the
input signals from $t_k$ to time $t_{k+j}$, is
computed by applying recursively, in the PAM and impulsive cases, by applying Equation (\ref{eqn-simp_evol}):
\begin{eqnarray}\label{eq:rec}
x_{k+j+1}&=&  A_{k+j+1,k} x_k
+\sum\limits_{i=k}^{k+j}   A_{k+j,i} 
B_{k,i}  U_{k,i},\qquad \end{eqnarray}
where the definition $A_{k,k}=\mathrm{I}$, $A_{k+1,k}=A_k$, and if $j> 0$ then $A_{k+j+1,k}=A_{k+j}A_{k+j-1}\hdots A_{k}$ has been used.
Define now $\mathcal{X}_k $ and $\mathcal{U}_k$ as a \emph{stack vector} of $N_p$ state and input vectors, respectively, spanning from time $t_{k+1}$ to time $t_{k+N_p}$ for the state and from time $t_{k}$ to time $t_{k+N_p-1}$ for the controls, where $N_p$ is the planning horizon:
\begin{eqnarray} \nonumber
\mathcal{X}_k &=& \left[ \begin{array}{c}x_{k+1}
\\  \vdots  \\
{x}_{k+N_p}
\end{array} \right],
\,
\mathcal{U}_k = \left[ \begin{array}{c}  U_k
\\  \vdots  \\
 U_{k+N_p-1}
\end{array} \right].
\,
\end{eqnarray}
Similarly, for the impulsive and PWM cases, define
\begin{eqnarray} \nonumber
{\Gamma}_k  &=& \left[ \begin{array}{c}\tau_{k}
\\  \vdots  \\
{\tau}_{k+N_p-1}
\end{array} \right],
\,{\Lambda}_k = \left[ \begin{array}{c}  \kappa_k
\\  \vdots  \\
 \kappa_{k+N_p-1}
\end{array} \right],\, \Upsilon_k=\left[ \begin{array}{c}  {\Gamma}_k
\\
{\Lambda}_k
\end{array} \right]
\end{eqnarray}
Then one can write
\begin{equation}
\label{eq:vep} \mathcal{X}_k=F_k x_k + G_k \mathcal{U}_k,
\end{equation}
where $G_k$  is a square, block lower triangular matrix of size $mN_p$,
defined as
\begin{equation}\label{G-def}
G_k  =  \left[ \begin{array}{cccc}  B_{k} & 0 & \cdots & 0 \\  A_{k+2,k+1} B_{k} &  B_{k+1}& \cdots & 0 \\
\vdots & \vdots & \ddots & \vdots \\
A_{k+N_p,k+1}B_{k} & A_{k+N_p,k+2} B_{k+1} & \hdots & B_{k+Np-1}
      \end{array} \right],
\end{equation}
this is, its non-null blocks are defined 
by $(G_k)_{jl}=  A_{k+j,k+l} B_{k+l-1}$, and the matrix $F_k$ is defined as:
\begin{equation}
F_k  =  \left[ \begin{array}{c}  A_{k+1,k} \\  A_{k+2,k} \\ \vdots \\   A_{k+N_p,k}
      \end{array} \right].
\end{equation}

It is important to note that, in the impulsive case, $G_k$ is a (nonlinear) function of $\Gamma_k$, whereas in the PWM case it is a (nonlinear) function of $\Upsilon_k$. To avoid lengthy expressions this dependence has been omitted. Another important remark is that, in the PWM case, $\mathcal U_k$ is fixed whereas in the other cases is the input variable

\section{Formulation of the planning problem}\label{sect-planning}

Next the planning problem is formulated, introducing the constraints and the objective function. The formulation is done for the three types of control signals.

\subsection{Constraints on the problem}\label{sec-constraints}
First constraints on the state and input are introduced. While only inequality are considered, equality constraints would be treated similarly.
\subsubsection{Inequality constraints on the state}\label{sec-state_constraints}

In this work it is assumed that the state is subject to inequality constraints along the planning horizon, which can vary as time advances. These can be formulated in general as $A_k \mathcal{X}_k    \leq b_k$, and using (\ref{eq:vep}), one reaches a expression in term of inputs, namely
\begin{equation}\label{eq:control_constraints}
A_k G_k \mathcal{U}_k  \leq b_k -A_k F_k x_k.
\end{equation}

%
%

\subsubsection{Input constraints}\label{sec-input_constraints}
Input constraints are different depending on the type of input.

In the PWM case, the inputs $\mathcal{U}_k$ are fixed, but the start time of impulse, $t_k+\tau_k$, and its end, $t_k+\tau_k+\kappa_k$, must be within the time interval both for negative and positive pulses. Thus
\begin{eqnarray}
0&\leq& \Gamma_k,\\
0&\leq& \Lambda_k,\\
\Gamma_k+ \Lambda_k&\leq & T,
\end{eqnarray}
which can be summarized as
\begin{eqnarray}
A^W \Upsilon_k \leq b^W.
\end{eqnarray}

In the PAM case, the inputs are limited above and below. Thus
\begin{equation}
\underline{U}_{PAM}\leq \mathcal{U}_k  \leq \overline{U}_{PAM}
\end{equation}

In the impulsive case, the inputs $u_k$ are limited above and below, but also the times of impulse, $t_k+\tau_k$, must be within the time interval. Thus
\begin{eqnarray}
\underline{U}_{IMP}&\leq& \mathcal{U}_k  \leq \overline{U}_{IMP},\\
0&\leq& \Gamma_k \leq T.
\end{eqnarray}

\subsection{Objective function}
The objective function to be minimized in the planning problem is a combination of the 1-norm of the control signal, which is denote das $J_\mathcal{U}$, (which gives an estimation of fuel consumption in case the control signal is thrust, see Section~\ref{sec:rendezvous}) and a weighted 2-norm of the state, which is denoted as $J_\mathcal{X}$, both taken over the planning horizon. Thus,
\begin{equation}
J_k=J_{\mathcal{U},k}+\alpha J_{\mathcal{X},k},
\end{equation}
where $\alpha$ is a positive constant that allows us to give a relative weight between input cost and state error. $J_{\mathcal{X},k}$ is computed as
\begin{equation}
J_{\mathcal{X},k}=\mathcal{X}_k^T  Q_k\mathcal{X}_k,
\end{equation}
for $Q_k>0$. Written in terms of the inputs and the starting point $x_k$, $J_{\mathcal{X},k}$ is
\begin{equation}
J_{\mathcal{X},k}=2x_k^T F_k^T Q_k G_k \mathcal{U}_k +\mathcal{U}_k^T G_k^T Q_k G_k \mathcal{U}_k ,
\end{equation}
an expression in which the constant term $x_k^T F_k^T Q_k F_k x_k$, which does not play a role in the planning optimization as it is constant for a given $x_k$, is neglected.

The value of $J_{\mathcal{U},k}$ does, however, depend on the control type.
\subsubsection{PWM control inputs}

For the case of PWM control inputs, using definition (\ref{eq:control_PWM_def}) it can be seen that the objective function $J_U(k)$ is given by:
\begin{eqnarray}
J_{\mathcal{U},k}&=& \sum_{j=k}^{k+N_p-1} \left[(U_j^+)^T \kappa_{j}^++(U_j^-)^T \kappa_{j}^-\right]
\nonumber \\ &=&
\mathcal U_k^T \Lambda_k
\nonumber \\ &=&
A^J_k \Upsilon_k,
\label{eq-costPWM}
\end{eqnarray}
with $A^J_k$ defined by blocks as
\begin{equation}\label{AJ-def}
A^J_k  =  \left[ \begin{array}{cc} 0 & 0 \\ 0 & U_k^T 
      \end{array} \right].
\end{equation}
The times $\Gamma_k$ where inputs start does not play a role in the cost function (only their duration $\Lambda_k$).

\subsubsection{PAM control inputs}

For the case of PAM control inputs, it can be seen that the objective function $J_U(k)$ is given by:
\begin{eqnarray}
J _{\mathcal{U},k}&=& \sum_{j=k}^{k+N_p-1}T \Vert U_j \Vert_1=T \Vert \mathcal{U}_k \Vert_1.\label{eq-costPAM}
\end{eqnarray}

\subsubsection{Impulsive control inputs}

For the case of impulsive control inputs, $J_U(k)$ is given by:
\begin{eqnarray}
J _{\mathcal{U},k}&=& \sum_{j=k}^{k+N_p-1} \Vert U_j \Vert_1=\Vert \mathcal{U}_k \Vert_1,\label{eq-cost}
\end{eqnarray}
where it should be noticed that, as in the PWM case, the location $\tau_{k,i}$ of the impulses does not play a role in the cost function.

\subsection{Planning optimization problem}
Now, for each of the input types, one can formulate a planning optimization problem starting from initial condition $x_k$ at time $t_k$, with a planning horizon of $N_p$, as follows.

\subsubsection{PWM control inputs}
For PWM control inputs, the planning optimization problem is formulated as

\begin{eqnarray} 
&\min_{\Upsilon_k} & 
2x_k^T F_k^T Q_k G_k(\Upsilon_k) \mathcal{U}_k  +\mathcal{U}_k^T G_k^T(\Upsilon_k ) Q_k G_k(\Upsilon_k ) \mathcal{U}_k
+\alpha A^J_k \Upsilon_k,\nonumber \\
\mbox{s. t.}
&  A_k G_k(\Upsilon_k) \mathcal{U}_k & \leq b_k -A_k F_k x_k.  \label{eq-controlcompPWM}\\
&A^W \Upsilon_k &\leq b^W.
\end{eqnarray}
Notice that $\mathcal{U}_k$ is known, and one has to compute the start and width of the pulses, contained in $\Upsilon_k$ (start and duration of pulses), which enter nonlinearly in the optimization problem. The dependence of $G_k$ on $\Upsilon_k$  has been made explicit.

\subsubsection{PAM control inputs}
For PAM control inputs, the planning optimization problem is formulated as
\begin{eqnarray} 
&\min_{\mathcal{U}_k} & 2x_k^T F_k^T Q_k G_k \mathcal{U}_k +\mathcal{U}_k^T G_k^T Q_k G_k \mathcal{U}_k
+\alpha T \Vert \mathcal{U}_k \Vert_1 , \nonumber \\
\mbox{s.t.}
&  A_k G_k \mathcal{U}_k & \leq b_k -A_k F_k x_k. \label{eq-controlcomp}\\
&\underline{U}_{PAM}&\leq \mathcal{U}_k \leq \overline{U}_{PAM}. \nonumber
\end{eqnarray}
\subsubsection{Impulsive control inputs}
For impulsive control inputs, the planning optimization problem is formulated as
\begin{eqnarray} \label{eq-controlcompimp}
&\min_{\mathcal{U}_k,\Gamma_k } & 
2x_k^T F_k^T Q_k G_k(\Gamma_k ) \mathcal{U}_k +\mathcal{U}_k^T G_k^T(\Gamma_k ) Q_k G_k(\Gamma_k ) \mathcal{U}_k
+\alpha  \Vert \mathcal{U}_k \Vert_1,\nonumber \\
\mbox{s. t.}
&  A_k G_k(\Gamma_k ) \mathcal{U}_k & \leq b_k -A_k F_k x_k.  \label{eq-controlcompimp}\\
&\underline{U}_{IMP}&\leq\mathcal{U}_k\leq \overline{U}_{IMP},\nonumber\\ \nonumber
&0&\leq \Gamma_k\leq T. \nonumber
\end{eqnarray}
The dependence of $G_k$ on $\Gamma_k$ (location of impulses) has been made explicit to emphasize that the optimization problem is nonlinear.

\section{PWM planning algorithm}\label{sect-method}
In this section the subindex $k$ is kept even though it does not play any role. For a ``pure'' planning problem, it could be set to zero. However, $k$ will be useful when defining the MPC algorithm in Section~\ref{sec-mpcrz}.

Consider now the problem (\ref{eq-controlcompPWM}), given $x_k$ and $\mathcal U_k$. Since the problem is nonlinear, one needs to design an algorithm to solve it. The planning algorithm is based on starting the problem using either the impulsive or PAM model. The algorithm is composed of the following steps.
\begin{description}
\item {\bf Step 1}. Solve either the PAM optimization problem (\ref{eq-controlcomp}), or the impulsive problem (\ref{eq-controlcompimp}) with a fixed $\Gamma_k$, to provide an initial guess of the PWM solution.
\item  {\bf Step 2}. The PAM or impulsive control inputs resulting from the optimization algorithm in Step 1 are converted to a sequence of PWM inputs, denote this initial sequence by $\Upsilon_k^0$. Set $i=0$.
 \item  {\bf Step 3}. The trajectory of the system with the PWM inputs $\Upsilon_k^i$ is computed analytically (if possible) or numerically by using equation (\ref{eq:vep}). Denote the trajectory by $\mathcal X_k^i$.
 \item  {\bf Step 4}. The system with PWM inputs is \emph{linearized} around $\mathcal X_k^i$, thus obtaining a linear, explicit plant with respect to \emph{increments}, denoted as $\Delta_k^i$, in the PWM inputs. Then a quadratic program can be posed and solved to find the increments that improve the cost function.
\item  {\bf Step 5}. The resulting solution $\Delta_k^i$ is used to improve the approximation towards the real solution, by setting $\Upsilon_k^{i+1}=\Upsilon_k^i+\Delta_k^i$. Increase $i$ by one and go back to Step 3. The process is iterated until the solution converges or time is up.
\end{description}

Next, all the steps in the algorithm are described.

\subsection{Step 1. Computation of PAM/impulsive control input}\label{sec-impulsi}

First, one has to choose if to find an initial guess using a PAM approach or an impulsive approach. The PAM guess is more suitable if one expects wide pulses, whereas the impulsive guess is best when the pulses are rather short.

If a PAM guess is chosen, it is computed from (\ref{eq-controlcomp}), setting $\overline{U}_{IMP}=T \mathcal{U}_k^+$ and $\underline{U}_{IMP}=T\mathcal{U}_k^-$, so that the solution can always be converted to PWM following the procedure of Section~\ref{sec-pwmfilt}. On the other hand, the impulsive guess is computed from (\ref{eq-controlcompimp}), setting $\overline{U}_{IMP}=\mathcal{U}_k^+$ and $\underline{U}_{IMP}=\mathcal{U}_k^-$. The impulsive guess also requires to set the impulse location $\Gamma_k$  to some pre-determined value, so only the impulse magnitude (which appears linearly in (\ref{eq-controlcompimp})) is unknown. Typical positions would be the middle of the interval (all entries of $\Gamma_k$ equal to $T/2$) or start of the interval ($\Gamma_k=0$).

\subsection{Step 2. Initial PWM solution: Adapting the PAM/impulsive solution}\label{sec-pwmfilt}
The PAM/impulsive solution from Section~\ref{sec-impulsi}, $\mathcal{U}$, is transformed to a PWM sequence of inputs, as follows:
%
%

\begin{enumerate}
\item From $\mathcal{U}$ extract $u^A_{j,i}$ (or $u^I_{j,i}$ if the initial solution is of impulsive type) for $j=k,\hdots,k+N_p-1$ and $i=1,\hdots,m$. Also extract $\tau_{j,i}$ if the initial solution is of impulsive type.
\item If the initial solution is of PAM type,  set
\begin{equation} 
\tau_{j,i}^+\hspace{-1pt}(t)=\left\{ \begin{array}{ll} \vspace{1pt}
 \frac{T u^A_{j,i}}{u^{W\pm}_{j,i}} , &  u^A_{j,i}>0,\\
0, &  u^A_{j,i}\leq0,
 \end{array} \right. \qquad \tau_{j,i}^- \hspace{-1pt}(t)=\left\{ \begin{array}{ll} \vspace{1pt}
 -\frac{T u^A_{j,i}}{u^{W\pm}_{j,i}} , &  u^A_{j,i}<0,\\
0, &  u^A_{j,i}\geq0,
 \end{array} \right.
\end{equation}
and if the initial solution is of impulsive type,
\begin{equation} 
\tau_{j,i}^+\hspace{-1pt}=\left\{ \begin{array}{ll} \vspace{1pt}
 \frac{u^I_{j,i}}{u^{W\pm}_{j,i}} , &  u^A_{j,i}>0,\\
0, &  u^I_{j,i}\leq0,
 \end{array} \right. \qquad \tau_{j,i}^- \hspace{-1pt}=\left\{ \begin{array}{ll} \vspace{1pt}
 -\frac{u^I_{j,i}}{u^{W\pm}_{j,i}} , &  u^A_{j,i}<0,\\
0, &  u^I_{j,i}\geq0,
 \end{array} \right.
\end{equation}
 
  \item In the PAM case, the PWM input should be centered in the interval: $\kappa_{j,i}^+=\frac{T- \tau_{j,i}^{+}}{2}$, $\kappa_{j,i}^-=\frac{T- \tau_{j,i}^{-}}{2}$. In the impulsive case, the PWM input should be centered around the chosen $\tau_{j,i}$ (corrected if necessary to avoid spillover), i.e. 
\begin{eqnarray} 
\kappa_{j,i}^+&=&\left\{ \begin{array}{ll} 
0, &  \tau_{j,i}-\frac{\tau_{j,i}^+}{2}<0,\\
T-\tau_{j,i}^+, & \tau_{j,i}+\frac{\tau_{j,i}^+}{2}>T,\\
\tau_{j,i}-\frac{\tau_{j,i}^+}{2} , &  \mathrm{otherwise},\\
 \end{array} \right. 
 \\
  \kappa_{j,i}^-&=&\left\{ \begin{array}{ll}
0, &  \tau_{j,i}-\frac{\tau_{j,i}^-}{2}<0,\\
T-\tau_{j,i}^-, & \tau_{j,i}+\frac{\tau_{j,i}^-}{2}>T,\\
\tau_{j,i}-\frac{\tau_{j,i}^-}{2} , &  \mathrm{otherwise},\\
 \end{array} \right.
 \end{eqnarray}
 \item From  $\tau_{j,i}^+$, $\tau_{j,i}^-$, $\kappa_{j,i}^+$, and $\kappa_{j,i}^-$, construct $\Gamma_k$ and $\Lambda_k$ and thus $\Upsilon_k$.
\end{enumerate}

The PWM signals $\Gamma_k,\Lambda_k$ constructed by this method produce a moderately similar (but not equal) output to the system driven by PAM or impulsive signals, but as time advances the output might considerably differ. See~\cite{Shieh1996673,Ieko1999123,Bernelli-Zazzera199864} for more details and other methods. In addition, the PWM results are not optimal (with respect to the PWM signals) and  they might not even verify the constraints. However, this solution is  only used as an initialization for the optimization algorithm proposed next. Denote as $\Upsilon_k^0$ the found solution and set $i=0$.

\subsection{Step 3. Computation of trajectories under PWM inputs}
For the current iteration $i$,  apply (\ref{eq:vep}) to compute the states of the system $\mathcal{X}_k^i$ at all times, with PWM inputs $\Upsilon_k^i$. The matrix $G_k$ might be needed to compute numerically if an explicit solution for the integrals (\ref{eqn-sol_pwm}) is not known or possible.

\subsection{Steps 4 and 5. Refined PWM solution: An optimization algorithm}\label{sec-refinement}
To linearize (\ref{eq:vep}) around inputs $\Upsilon_k^i$,
notice from (\ref{eqn_BWdef}) that
\begin{eqnarray}
\frac{\partial}{\partial \tau_{k,i}} B^W_{k,i}(\tau_{k,i},\kappa_{k,i})&=& \frac{\partial}{\partial \tau_{k,i}} \int_{t_k+\tau_{k,i}}^{t_k+\tau_{k,i}+\kappa_{k,i}} \Phi(t_{k+1},s) B_{i}(s) ds,
\nonumber \\Ê
&=&\Phi(t_{k+1},t_k+\tau_{k,i}+\kappa_i) B_i(t_k+\tau_{k,i}+\kappa_{k,i})  
\nonumber \\
&&-\Phi(t_{k+1},t_k+\tau_{k,i}) B_i(t_k+\tau_{k,i}) ,
\end{eqnarray}
and
\begin{eqnarray}
\frac{\partial}{\partial \kappa_{k,i}} B^W_{k,i}(\tau_{k,i},\kappa_{k,i})&=& \frac{\partial}{\partial \kappa_{k,i}} \int_{t_k+\tau_{k,i}}^{t_k+\tau_{k,i}+\kappa_{k,i}} \Phi(t_{k+1},s) B_{i}(s) ds,
\nonumber \\Ê
&=&\Phi(t_{k+1},t_k+\tau_{k,i}+\kappa_i) B_i(t_k+\tau_{k,i}+\kappa_{k,i}) .
\end{eqnarray}
Thus, (\ref{eqn-simp_evol}) can be explicitly linearized around some given $\bar \tau_{k}^+$, $\bar \tau_{k}^-$ and $\bar \kappa_{k}^+$, $\bar \kappa_{k}^-$, reaching
\begin{equation}
x_{k+1}=A_k x_k+B_k U_{k}+B^{\Delta \tau}_k \delta \tau_k+B^{\Delta \kappa}_k \delta \kappa_k,\label{eqn-simp_evol_lin}
\end{equation}
where $B_k$ is computed with $\bar \tau_{k}^+$, $\bar \tau_{k}^-$, $\bar \kappa_{k}^+$, and $\bar \kappa_{k}^-$, and  define
\begin{eqnarray}
B^{\Delta \tau}_k&=& \left[ \begin{array}{cc} B_{k}^{\delta \tau^+} & B_{k}^{\delta\tau^-} \end{array} \right],\quad B^{\Delta \kappa}_k= \left[ \begin{array}{cc} B_{k}^{\delta\kappa^+} & B_{k}^{\delta\kappa^+}
\end{array} \right],\\
{\delta \tau}_k &=& \left[ \begin{array}{c} \bar \tau_{k}^+- \tau_{k}^+ \\ \bar \tau_{k}^- -\tau_{k}^- 
\end{array} \right],\quad {\delta \kappa}_k = \left[ \begin{array}{c} \bar \kappa_{k}^+- \kappa_{k}^+ \\ \bar \kappa_{k}^-- \kappa_{k}^-
\end{array} \right], \label{eqn-bkdeltas}
\end{eqnarray}
where the $i$-th entries of the $B_k^\delta$ matrices in (\ref{eqn-bkdeltas}) are given, respectively, by
\begin{eqnarray}
(B_{k}^{\delta \bar \tau^+} )_i&=&\Phi(t_k+T,t_k+\bar \tau_{k,i}^++\bar \kappa_{k,i}^+) B_{i}(t_k+\bar \tau_{k,i}^++\bar \kappa_{k,i}^+)  u^{W+}_{k,i} \nonumber \\ && - \Phi(t_k+T,t_k+\bar \tau_{k,i}^+) B_{i}(t_k+\bar \tau_{k,i}^+)  u^{W+}_{k,i}  \\
(B_{k}^{\delta \bar \tau^-} )_{i}&=&-\Phi(t_k+T,t_k+\bar \tau_{k,i}^-+\bar \kappa_{k,i}^-) B_{i}(t_k+\bar \tau_{k,i}^-+\bar \kappa_{k,i}^-)  u^{W-}_{k,i}
\nonumber \\ &&
+\Phi(t_k+T,t_k+\bar \tau_{k,i}^-) B_{i}(t_k+\bar \tau_{k,i}^-)  u^{W-}_{k,i}\\
(B_{k}^{\delta \bar \kappa^+} )_{i}&=&\Phi(t_k+T,t_k+\bar \tau_{k,i}^++\bar \kappa_{k,i}^+) B_{i}(t_k+\bar \tau_{k,i}^++\bar \kappa_{k,i}^+)  u^{W+}_{k,i},\\
(B_{k}^{\delta \bar \kappa^-} )_{i} &=&-\Phi(t_k+T,t_k+\bar \tau_{k,i}^-+\bar \kappa_{k,i}^-) B_{i}(t_k+\bar \tau_{k,i}^-+\bar \kappa_{k,i}^-)  u^{W-}_{k,i}
\end{eqnarray}
Thus, defining stack vectors with the increments in the PWM variables at step $i$ as
\begin{eqnarray} \nonumber
{\Delta \Gamma}_k^i  &=& \left[ \begin{array}{c}\delta \tau_{k}
\\  \vdots  \\
{\delta} \tau_{k+N_p-1}
\end{array} \right],\quad {\Delta \Lambda}_k^i  = \left[ \begin{array}{c}\delta \kappa_{k}
\\  \vdots  \\
{\delta} \kappa_{k+N_p-1}
\end{array} \right],
\end{eqnarray}
and grouping all increments as
\begin{eqnarray} \nonumber
{\Delta }_k^i  &=& \left[ \begin{array}{c}{\Delta \Gamma}_k^i  \\
{\Delta \Lambda}_k^i  
\end{array} \right],\quad B_k^\Delta  = \left[ \begin{array}{c}B^{\Delta \tau}_k \\ B^{\Delta \kappa}_k
\end{array} \right],
\end{eqnarray}
and defining $G_k^{\Delta}$ as in (\ref{G-def}), i.e.,
\begin{equation}\label{G-def}
G^{\Delta }_k  =  \left[ \begin{array}{cccc}  B^\Delta_{k} & 0 & \cdots & 0 \\  A_{k+2,k+1} B^\Delta_{k} &  B^\Delta_{k+1}& \cdots & 0 \\
\vdots & \vdots & \ddots & \vdots \\
A_{k+N_p,k+1}B^\Delta_{k} & A_{k+N_p,k+2} B^\Delta_{k+1} & \hdots & B^\Delta_{k+Np-1}
      \end{array} \right],
\end{equation}
one can write
\begin{equation}
\label{eq:vep_PWM} \mathcal{X}^i_k\approx F_k x_k + G_k(\Upsilon_k^i) \mathcal{U}_k+G_k^{\Delta}(\Upsilon_k^i) \Delta_k^i,
\end{equation}
The state
constraints (\ref{eq:control_constraints}) 
become
\begin{equation}\label{eq:control_constraints_PWM2}
A_k G_k^{\Delta}(\Upsilon_k^i) \Delta_k^i  \leq b_k -A_k F_k x_k-G_k(\Upsilon_k^i) \mathcal{U}_k.
\end{equation}
The constraints on $\Delta \Gamma_k^i$ and $\Delta \Lambda_k^i$ are as follows:
\begin{eqnarray}
-A^W \Delta_k^i& \leq& b^W-A^W \Upsilon_k^i,\\
\Delta_{max} &\leq & \Delta_k^i \leq \Delta_{max},\label{eqn:deltamax}
\end{eqnarray}
where the last constraint (\ref{eqn:deltamax}) is used to avoid large variations that might make the linearization approximation to fail. All these constraints are might be summarized as
\begin{equation}
A^\Delta  \Delta_k^i\leq b^\Delta
\end{equation}
Finally, the objective function can be rewritten in terms of $\Delta_k^i$ as 
$J_k^i=J_k^i(\Upsilon_k^i)+J_k^\Delta(\Upsilon_k^i,\Delta_k^i)$, where
\begin{eqnarray}
J_k^i&=&2x_k^T F_k^T Q_k G_k(\Upsilon_k^i) \mathcal{U}_k  +\mathcal{U}_k^T G_k^T(\Upsilon_k^i ) Q_k G_k(\Upsilon_k^i ) \mathcal{U}_k
+\alpha A^J_k \Upsilon_k^i,\label{eq-costPWM1}\\
J_k^\Delta&=&2(x_k^T F_k^T+ \mathcal{U}_k^TG_k^T(\Upsilon_k^i)) Q_k G_k^\Delta(\Upsilon_k^i) \Delta_k^i  +(\Delta_k^i)^T G_k^{\Delta T}(\Upsilon_k^i ) Q_k G_k^\Delta(\Upsilon_k^i )\Delta_k^i
\nonumber \\ &&
+\alpha A^J_k \Delta_k^i.\label{eq-costPWM1}
\end{eqnarray}
Noting that $J_k^\Delta$ is quadratic in $\Delta_k^i$, a quadratic optimization program with linear restriction, formulated on the output increments, can be posed as follows:
\begin{eqnarray} \label{eq-controlcomprob_PWM}
\min_{\Delta_k^i} &J_k^\Delta(\Upsilon_k^i,\Delta_k^i) \\
\mbox{s. t.:}
& A_k G_k^{\Delta}(\Upsilon_k^i) \Delta_k^i & \leq b_k -A_k F_k x_k-G_k(\Upsilon_k^i) \mathcal{U}_k. , \nonumber \\
&A^\Delta  \Delta_k^i &\leq b^\Delta.\nonumber
\end{eqnarray}
The solution $\Delta_k$ is used to recompute new PWM inputs, 
$\Upsilon^{i+1}_k=\Upsilon^{i}_k+\Delta_k^i$.
Then the linearization process can be repeated around the new $\Upsilon^{i+1}_k$, refining the solution in each iteration.

\section{Model Predictive Control with PWM inputs}\label{sec-mpcrz}
In this section, building upon the trajectory planning algorithm of Section~\ref{sect-method}, which is open-loop and has a finite time-horizon, a closed-loop algorithm is developed based on the ideas of model predictive control (also known as receding horizon control). Model predictive control closes the loop by simply re-planning the maneuver at each time step, after applying just the set of control inputs corresponding to the first time step, and keeps looking ahead $N_p$ time steps. Thus, the algorithm starts at $k=0$ and is repeated for each $k$. The re-planning is done from the actual position at each time step, which seldom coincides with the planned position due to disturbances, thus effectively closing the loop.

However, except at the start, it is not necessary to repeat all the steps of Section~\ref{sect-method}. Since the new position should be close to the planned one, one can apply the linearization scheme of the planning algorithm starting from the last available linearization. The MPC algorithm is summarized next:
\begin{description}
\item {\bf Step 1}. At time step $k=0$ and starting from $x_0$ apply the Planning algorithm of Section~\ref{sect-method}, obtaining a set of PWM inputs $\Upsilon_0$ that would optimize the planning problem (\ref{eq-controlcompPWM}) for the next $N_p$ time steps, if there were no disturbances.
\item  {\bf Step 2}. Apply impulses corresponding to the first time instant; save the rest of impulses. Set $k=1$
 \item  {\bf Step 3}. One arrives at $x_k$, which probably is not the intended value of the state at time $k$ but close. Thus re-planning is necessary.
 \item  {\bf Step 4}. For re-planning,  apply the planning algorithm of Section~\ref{sect-method}. However, to avoid the initial step of having to use a PAM or impulsive model and compute an initial guess, use instead as an initial guess the impulses of $\Upsilon_{k-1}$ that were not used (all of them except those corresponding to time $k-1$) and guess the remaining impulses (at the end) as zeros. In this way, form an initial guess $\Upsilon^0_k$.
 \item  {\bf Step 4}. Apply the linearization algorithm of Section~\ref{sec-refinement} using $\Upsilon^0_k$ as initial guess to obtain, after iterating, a new set of impulses $\Upsilon_k$. Apply the set of impulses corresponding to time $k$. Save the rest of impulses.
\item  {\bf Step 5}. Repeat step 3.
\end{description}

\section{Example application: Spacecraft Rendezvous}\label{sec:rendezvous}
Rendezvous of spacecraft is the controlled close encounter of two (or more) space vehicles. This work assumes just two vehicles, one of which is the \emph{target vehicle} (which is in a known orbit, and considered passive) and the other is the \emph{chaser spacecraft}, which begins from a known position and maneuvers until very close to target. Only close range rendezvous~\cite{wigbert} is considered, which starts at hundreds of meters and ends when the chaser is very close to  target (a few meters with speeds of centimeters per second).

There are numerous mathematical models for
spacecraft rendezvous; which one  should be used depends on the
parameters of the scenario. In~\cite{Carter} a survey of numerous
mathematical models for spacecraft rendezvous can be found.

For instance, if the target is orbiting in a \emph{circular}
Keplerian orbit, the general equations of the relative movement between an active chaser spacecraft close to a passive target vehicle are linear time-invariant Hill-Clohessy-Wiltshire (HCW) equations (introduced in~\cite{Hill} and~\cite{CW}). While these equations are frequently used in the literature, it must be  noted that, in many situations, the HCW equations are
not accurate. For instance, if the target vehicle is moving in a Keplerian
\emph{eccentric} orbit (see~\cite{How2002}) or if some orbital
perturbations are taken into account (see for
example~\cite{Humi}). A more complex model, the Tschauner-Hempel model (see~\cite{T-H} or~\cite{Carter}) assumes that the target vehicle is passive and
moving along an elliptical orbit with semi-major axis $a$ and eccentricity $e$. The system equations are linear time-varying and cannot be exactly integrated in time to obtain a discrete transition model; however, if one substitutes the time $t$ by the \emph{eccentric} anomaly of the target orbit, $E$, it is possible to obtain explicit expressions for the system evolution in the PWM, impulsive, and PAM  actuation cases. This will be the model considered in this work. The model can be expressed in cartesian coordinates, but also in the so-called relative orbital elements (see, e.g.,~\cite{gaias2} or~\cite{sinclair2014}). The former has been chosen for simplicity.

Let us first establish some notation. Define the orbital mean motion  $n=\sqrt{\frac{\mu}{a^3}}$, where $\mu$ is the gravitational parameter of the central body around which the target spacecraft is orbiting.

Now, note that $t$ and $E$ are related in a one-to-one fashion by using Kepler's equation:
\begin{equation}
n(t-t_p)=E-e\sin E,
\end{equation}
where $t_p$ is the time at periapsis used as a starting point to measure the eccentric anomaly $E$. The time $t_p$ is chosen such that it is equal or less than the starting time which is denoted as $t_0$ (subtracting, if necessary, any number of orbital periods). Kepler's equation is not analytically invertible, but its inverse can be found numerically with any desired degree of precision
 (see any Orbital Mechanics reference, such as~\cite{Wie:1998}). Denote its inverse by the function $K$, i.e. $E=K(t)$. Denote by $E_0$ the true anomaly corresponding to $t_0$, this is, $E_0=K(t_0)$. Then, $E_k=K(t_k)=K(t_0+kT)$, where $T$ is the sampling time (not to be confused with the orbital period). Call as $r_{x,k}$, $r_{y,k}$, and $r_{z,k}$ the position of the chaser in a local--vertical/local--horizontal (LVLH) frame of reference fixed on the center of gravity of the target vehicle at time $t_k$. In the (elliptical) LVLH frame, $x$ refers to the radial position, $z$ to the out-of-plane position (in the direction of the orbital angular momentum), and $y$ is perpendicular to these coordinates (not necessarily aligned with the target velocity given that its orbit is not circular).  The velocity and inputs of the chaser in the LVLH frame at time $t_k$ are denoted, respectively, by $v_{x,k}$, $v_{y,k}$, and $v_{z,k}$, and by $u_{x,k}$, $u_{y,k}$, and $u_{z,k}$.

\begin{figure}
\begin{center}
\includegraphics[width=0.8\columnwidth]{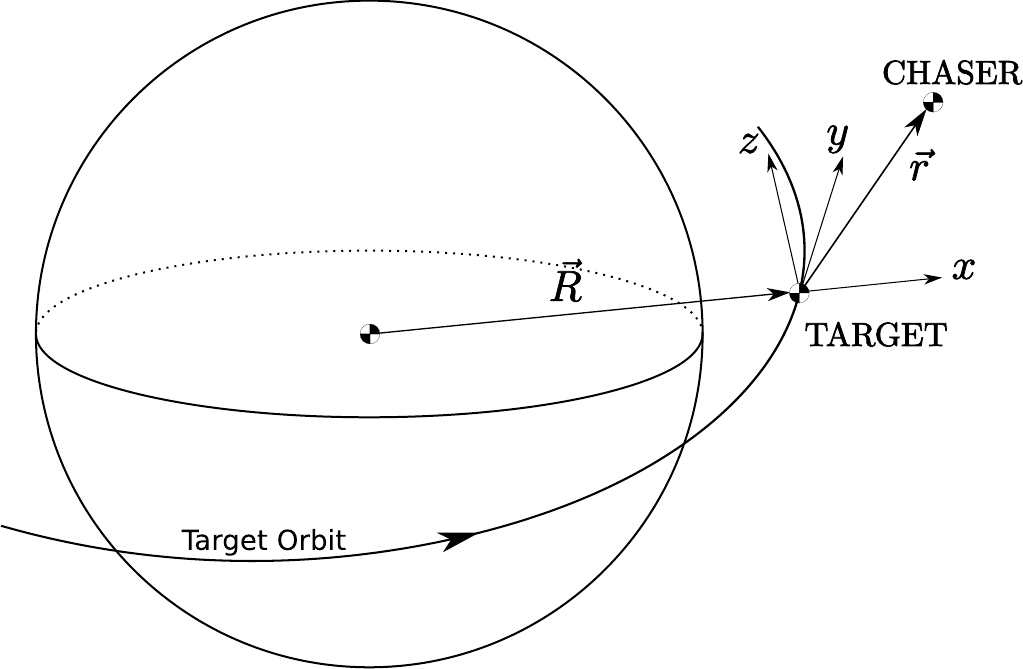}    
\caption{LVLH frame.}
\label{fig-lvlh frame}
\end{center}
\end{figure}

If there is no actuation (i.e. $u_{x,k}=u_{y,k}=u_{z,k}=0$), the resulting transition equation was obtained in a simple form in~\cite{Y-A} as follows:
\begin{eqnarray}
 x_{k+1}=\Phi(t_{k+1},t_{k})x_k
\end{eqnarray}
where
\begin{eqnarray}
x_k&=&\left[r_{x,k}~r_{y,k}~r_{z,k}~v_{x,k}~v_{y,k}~v_{z,k}\right]^T\hspace{-2pt},
\end{eqnarray}
and where 
\begin{equation}
\Phi(t_{k+1},t_{k})=Y_{K(t_{k+1})}Y^{-1}_{K(t_{k})},\label{eqn-fundamental}
\end{equation}
with $Y_{t_k}$ being the fundamental matrix solution of the Tschauner-Hempel model, which are expressed in~\cite{Y-A} as a function of \emph{true} anomaly $\theta$. However there is a one-to-one relation between $E$ and $\theta$ given by
\begin{equation}
\tan \frac{\theta}{2}=\sqrt{\frac{1+e}{1-e}} \tan \frac{E}{2},
\end{equation}
 which is exploited in the sequel. The explicit expression of the matrices\footnote{These expressions slightly differ from~\cite{Y-A} because the two transformation matrices that appear in that paper have been pre-multiplied; also, the reference axes are not the same.} is found in (\ref{eqn-Y}) and (\ref{eqn-Yinv}),
\begin{eqnarray}
\hspace{-15pt}&\hspace{-5pt}&\hspace{-5pt}Y_E=\left[\scriptsize \begin{array}{cccccc}
s & 0 & 0 & 2/\rho-3esJ &-c & 0\\
   c(1+1/\rho) & 1/\rho& 0& -3\rho J & s(1+1/\rho)& 0\\
   0 &0 &c/\rho &0 &0 &s/\rho\\
   \alpha  \rho^2 c& 0& 0& \alpha (-es-3e\rho^2Jc) & \alpha \rho^2s& 0\\
   \alpha  s(-1-\rho^2)& \alpha es& 0 &\alpha \rho(3es\rho J-3)& \alpha (c+e+c\rho^2) &0\\
   0 &0 &-s\alpha  &0 &0 &(c+e)\alpha  \end{array} \right],\label{eqn-Y}\normalsize
   \end{eqnarray} 
   \begin{eqnarray}
\hspace{-15pt}  &\hspace{-5pt}&\hspace{-5pt}Y^{-1}_E=\frac{3J}{(1-e^2)} \left[\scriptsize \begin{array}{cccccc}
  e \rho^2(1+\rho) & -e^2\rho^2 s & 0 & e^2s/\alpha & e\rho/\alpha& 0\\
   \rho^2(1+\rho) & -e\rho^2 s& 0& es/\alpha & \rho/\alpha& 0\\
   0 &0 &0 &0 &0 &0\\
   0 & 0 & 0 & 0 & 0& 0\\
   0 & 0 & 0 & 0 & 0& 0\\
   0 & 0 & 0 & 0 & 0& 0  \end{array} \right]+\frac{1}{(1-e^2)}
  \nonumber \\ \hspace{-15pt}   &\hspace{-5pt}&\hspace{-5pt}\times \hspace{-5pt}
 \left[\scriptsize\begin{array}{cccccc}
-s(\rho^2+2\rho+e^2) & es^2(1+\rho) & 0& \frac{c-2e/\rho}{\alpha } & -\frac{s(\rho+1)}{\rho\alpha } & 0\\
   -es(1+\rho)^2 & \rho^2(1-ce)+e^2s^2 & 0 & \frac{ec-2/\rho}{\alpha } & -\frac{es(\rho+1)}{\rho\alpha }& 0\\
   0 & 0 & (c+e)(1-e^2) & 0 & 0 & \frac{-s(1-e^2)}{\alpha \rho}\\ \vspace{1pt}
   \rho^2(1+\rho) & -es\rho^2 & 0 & \frac{es}{\alpha } & \frac{\rho}{\alpha } & 0\\
   3\rho(c+e)-e\rho s^2 & -esc(1+\rho)-e^2s & 0 & \frac{s}{\alpha } & \frac{c(\rho+1)+e}{\alpha \rho} & 0\\
    0 & 0 & s(1-e^2) & 0 & 0 & \frac{c(1-e^2)}{\alpha \rho}
     \end{array} \right]\quad
     \label{eqn-Yinv} 
\end{eqnarray} 
where the following symbols are used (expressed in terms of $E$):
\begin{eqnarray}
\rho&\hspace{-1pt}=\hspace{-1pt}&\frac{1-e^2}{1-e\cos E},\, s=\frac{\sqrt{1-e^2}\sin E}{1-e\cos E},c=\frac{\cos E-e}{1-e\cos E},\\
J&\hspace{-1pt}=\hspace{-1pt}&\alpha\frac{E-\hat E-e(\sin E-\sin \hat E)}{(1-e^2)^{3/2}},\hspace{-1pt}\alpha \hspace{-1pt}=\hspace{-1pt}\frac{n}{(1-e^2)^{3/2}}\label{eqn-J},\quad
\end{eqnarray}
where $\hat E$ in (\ref{eqn-J}) can be substituted by zero or any other desired reference value of $E$. For instance, if when evaluating (\ref{eqn-fundamental}) one chooses $\hat E=E_{k}=K(t_{k})$, then for $Y^{-1}_{K(t_{k})}$ one gets $J=0$ and the first matrix in (\ref{eqn-Yinv}) becomes zero.

Using  (\ref{eqn-fundamental}), one gets $A_k$ in (\ref{eqn-simp_evol}) explicitly, as well as $B_k$ in the impulsive case (explicitly defined in in terms of $\Phi(t_{k+1},t_{k})$). To obtain the $B_k$ matrix in the PWM and PAM cases, one needs to solve (\ref{eqn_BWdef}) or (\ref{eqn-defBA}), respectively, which involves an integral. Defining
\begin{equation}
b_i(r_1,r_2,r_3)= \int_{r_1}^{r_2} \Phi(r_3,s) B_{i}(s) ds,
\end{equation}
one has that, from (\ref{eqn_BWdef}),
\begin{equation}
B^W_{k,i}(\tau_{k,i},\kappa_{k,i})= b_i(t_k+\tau_{k,i},t_k+\tau_{k,i}+\kappa_{k,i},t_{k+1}),\end{equation}
and, from (\ref{eqn-defBA}), 
\begin{equation}
B^A_{k,i}= b_i(t_k,t_{k+1},t_{k+1}).\end{equation}
To compute the $b_i$'s, the following integral is needed
\begin{equation}
 \int  Y^{-1}_{K(t)}C_{i+3} dt, \label{eqn-inh}
\end{equation}
where $C_i$ is a 6-element column vector of zeros with a value of one at row $i$.  For the computation, define the functions $f_i(t)$, for $i=1,2,3$, as the indefinite integrals of (\ref{eqn-inh}), in terms of eccentric anomaly
\begin{equation}
f_i(E)= \int  Y^{-1}_{E}C_{i+3} \frac{1-e\cos E}{n}dE. \label{eqn-inh2}
\end{equation}
Once the $f_i$'s are computed, one finds the $b_i$'s as
\begin{equation}
b_i(r_1,r_2,r_3)=Y_{K(r_3)} \left(f_i(K(r_2))-f_i(K(r_1))\right)
\end{equation}
Inserting the expression of (\ref{eqn-Yinv}) in (\ref{eqn-inh2}) and integrating, one obtains
\begin{eqnarray}
 f_1\hspace{-4pt}&=&\hspace{-4pt}\frac{\left(1-e^2\right)^{-7/2}}{2\alpha^2}
\left[\scriptsize
\begin{array}{c}
2(1+6e^2)S-3e(2+e^2)E+e^2Ch+e^3\hat S/2 \\
2e(8-e^2)S-(4+7e^2-2e^4)E+ eCh+e^2\hat S/2 \\
  0\\
    -2e C    \left(1-e^2\right)^{3/2}\\
  -2C   \left(1-e^2\right)^{3/2}\\
  0
  \end{array}
\right]\hspace{-3pt},\,\quad\\
  f_2\hspace{-4pt}&=&\hspace{-4pt}\frac{\left(1-e^2\right)^{-3}}{2 \alpha^2} \left[\scriptsize
\begin{array}{c}
C(4(1+e^2)-eC)-eEh-3eE^2\\
eC(10-2e^2-eC)-Eh-3E^2\\
 \\
 0
 \\
2E \left(1-e^2\right)^{3/2} \\
    \sqrt{1-e^2}  (4S-e(3 E+\hat S/2 ) )
 \\ 0 \end{array}
\right]\hspace{-3pt},\quad\\
f_3\hspace{-4pt}&=&\hspace{-4pt}\frac{\left(1-e^2\right)^{-5/2}}{4 \alpha^2} \left[\scriptsize
\begin{array}{c} 0 \\ 0 \\
  2 \sqrt{1-e^2} C  (2-e C )\\
  0 \\ 
  0 \\
4 \left(e^2+1\right) S -e (6 E+\hat S)
\end{array}
\right],
\end{eqnarray}
where $S=\sin (E)$, $\hat S=\sin (2E)$, $C=\cos E$, $h=6\alpha(\hat E-E-e\sin (\hat E))$. Similar expressions for the $B$ matrices can be found in~\cite{ankersen-book}, however using a slightly different definition of reference axes.

Note that, using these formulas, it is possible to express (\ref{eqn-simp_evol}) explicitly for all actuation types. This greatly speeds up the algorithms.

\subsection{Constraints for the rendezvous problem}
\begin{figure}[h]
\centering
\includegraphics[width=0.8\columnwidth]{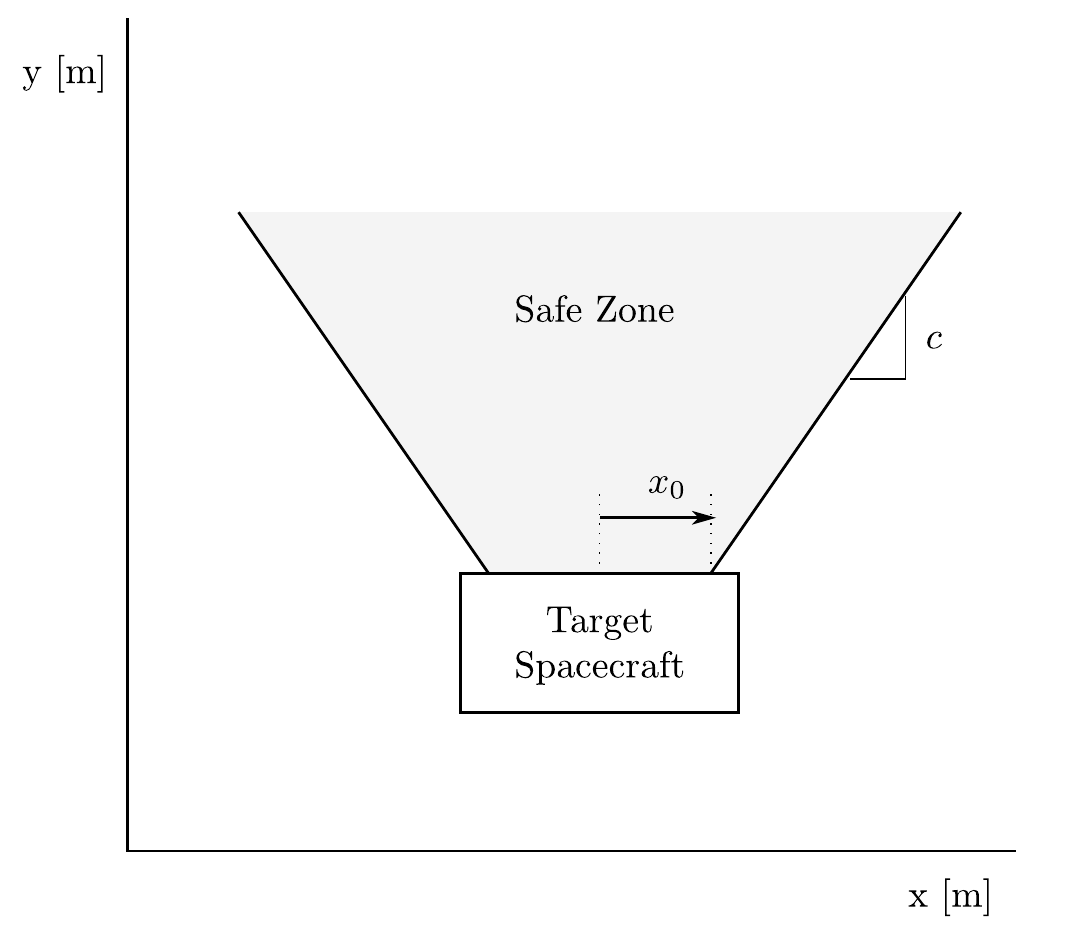}    
\centering
\caption{Line of Sight region.}
\label{fig-los}
\end{figure}

Besides the input constraints (which were given in Section~\ref{sec-input_constraints}), the inequality state constraints which were generically specified in Section~\ref{sec-state_constraints} are, in general, related to safety and sensing considerations (see e.g.~\cite{How2008}). In this work, it is considered that during rendezvous the chaser vehicle has to remain inside a line of sight (LOS) area. To simplify the constraint\footnote{More complicated constraints could be considered, see~\cite{gavilan2012} for examples including a rotating LOS constraint.}, in this work a 2-D LOS area is used as shown in Figure~\ref{fig-los}. This LOS region is the intersection of a cone, given by the equations $r_y\geq c_{LOS}(r_x-r_{x_0})$ and $r_y \geq -c_{LOS}(r_x+r_{x_0})$, and the region $r_y\geq0$.

The LOS constraint is $A_{LOS} x_k   \leq b_{LOS}$, where
\begin{equation}
A_{LOS} =
\left[ \begin{array}{cccccc}
        0 & -1 & 0 & 0 & 0 & 0 \\
        c_{LOS} & -1& 0 & 0 & 0 &  0 \\
        -c_{LOS} & -1 & 0 & 0 & 0 &  0
       \end{array}\right],\quad b_{LOS}=\left[ \begin{array}{c}
0 \\ c_{LOS}r_{x_0} \\ c_{LOS}r_{x_0} \end{array} \right]. \label{eq:state_constrains}
 \end{equation}

Using the compact formulation that was developed in Section~\ref{sect-model}, the constraints equations for
the state can be rewritten as:
\begin{equation}
\label{eq:state_constraints_comp}
A_c \mathcal{X} \leq b_c,
\end{equation}
where $A_c$ and $b_c$ are given by:
\begin{eqnarray}
A_c &=& \left[ \begin{array}{cccc}
A & & & \\
 &A & & \\
 & & \ddots & \\
 & & & A\\
\end{array}
 \right],\,
b_c = \left[ \begin{array}{c}
 b_{LOS} \\
  b_{LOS} \\
 \vdots \\
  b_{LOS}\\
\end{array}
 \right].
\end{eqnarray}

Then, using equation (\ref{eq:vep}), one can reformulate  the LOS
constraints as constraints for the control signals, starting at time step $t_k$, in the
following way:
\begin{equation}\label{eq:control_constraints}
A_c G_k \mathcal{U}_k   \leq b_c -A_c  F_k x_0.
\end{equation}
\section{Simulation Results} \label{sect-simu}

\begin{figure}[!t]
\begin{center}
\includegraphics[width=0.9\columnwidth]{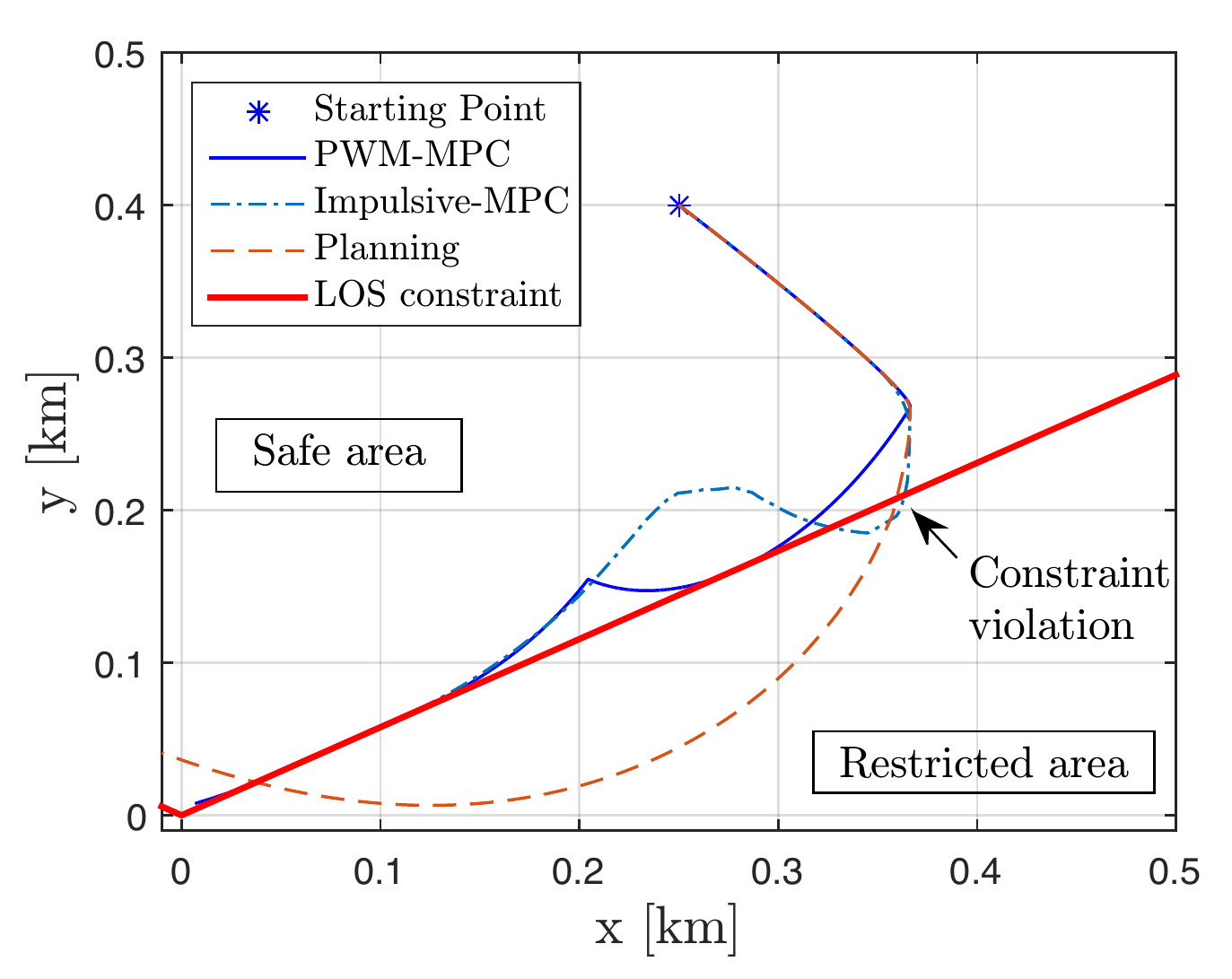}    

\caption{System trajectories in the target orbital plane: open-loop PWM inputs computed from impulsive solution (dashed), closed-loop Model Predictive Control with PWM inputs using impulsive model (dot-dashed), and closed-loop Model Predictive Control with PWM inputs using the PWM planning algorithm (solid).}
\label{fig-traj}
\end{center}
\end{figure}
\begin{figure}[!t]
\begin{center}
\includegraphics[width=0.9\columnwidth]{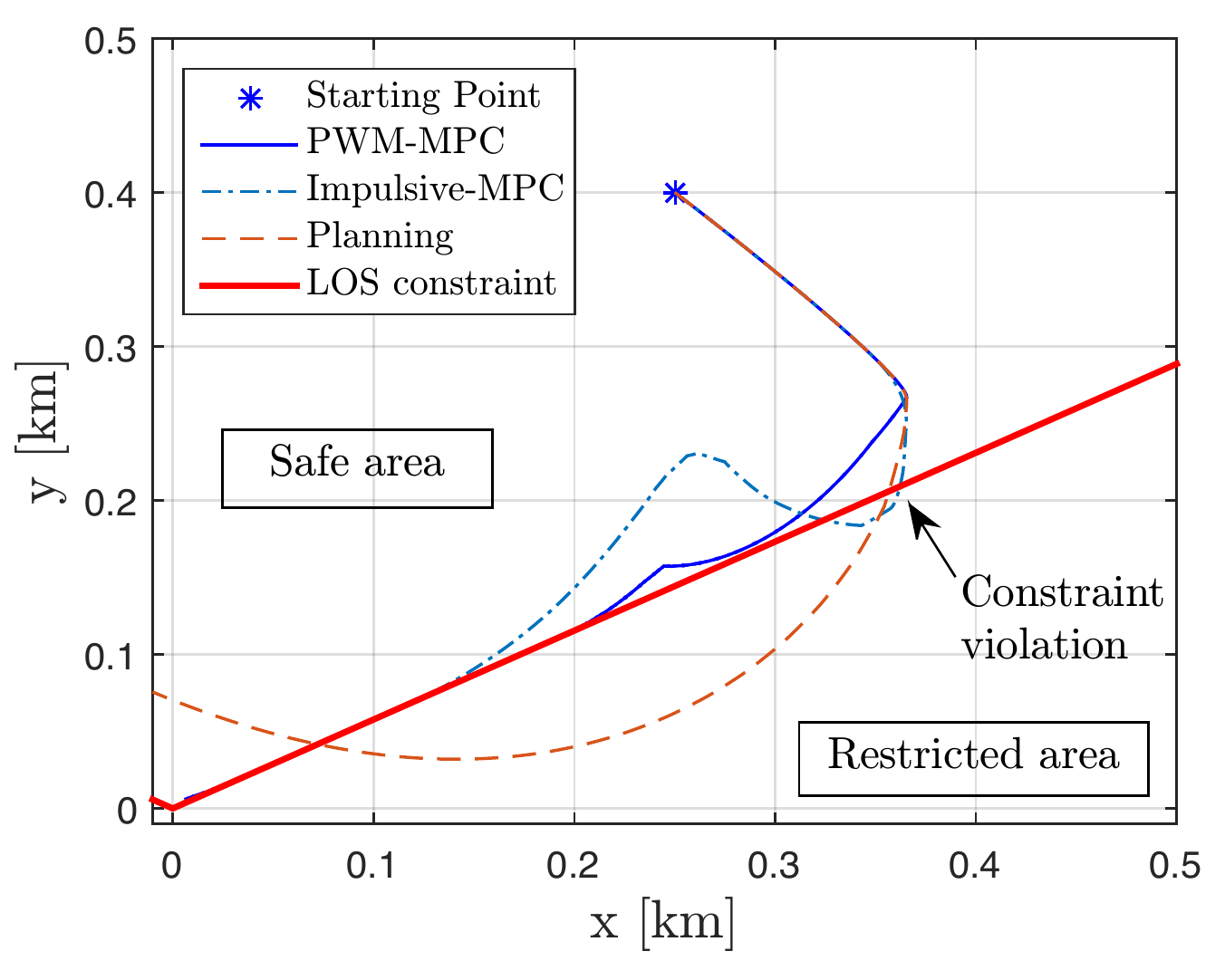}    

\caption{System trajectories in the target orbital plane, with inexact orbit model: open-loop PWM inputs computed with the planning algorithm (dashed), closed-loop Model Predictive Control with PWM inputs using impulsive model (dot-dashed), and closed-loop Model Predictive Control with PWM inputs using the PWM planning algorithm (solid).}
\label{fig-dist}
\end{center}
\end{figure}

For simulations the following values have been used: $N_p=50$ as planning horizon, $T=60~\mathrm{s}$, and $\bar u=10^{-1}~\mathrm{N/kg}$. The target orbit has $e=0.7$ and perigee altitude $h_p=500~\mathrm{km}$. Initial conditions were $\theta_0=45^\mathrm{o}$,
${\bf r}_0=[0.25\,\,0.4\,\,-0.2]^T~\mathrm{km},\,{\bf v}_0=[0.005\,\,-0.005\,\,-0.005]^T~\mathrm{km/s}$. The LOS constraint  (see~\cite{IFAC_PWM}) is defined by $x_0=0.001~\mathrm{km}$ and $C_{LOS}=\tan 30^\mathrm{o}$. For the cost function, $\alpha$ has been set to $10^3$ and $Q_k$ as
\begin{equation}
Q_k = \left[\begin{array}{ccc}
R_{k+1} & & \\
& \ddots & \\ &&
R_{k+N_p}
\end{array}\right],
\end{equation}
where $R_k$ is defined as
\begin{eqnarray} \label{eq-defR}
R_k= h (k-k_a)
 \left[
 \begin{array}{c c }
 \mathrm{Id}_{3\times 3} & {\Theta}_{3\times 3} \\
 {\Theta}_{3\times 3}  & {\Theta}_{3\times 3} \\
 \end{array}
 \right].
%
\end{eqnarray}
In (\ref{eq-defR}), $h$ is the step function, $k_a$ is the desired
arrival  time for  rendezvous, and $
\mathrm{Id}_{3\times 3}$, ${\Theta}_{3\times 3}$ are respectively
the identity matrix and a matrix full of zeros, both of order 3
by 3. The reason for choosing (\ref{eq-defR}) is that it is desired to arrive
at the origin at time $k_a$ (and \emph{remain} there) and at the
same time minimize the control effort.

In the simulations three algorithms were considered: first, an impulsive open-loop trajectory planner, as described in Section~\ref{sec-impulsi}. Next, closed-loop simulations using MPC but considering impulsive instead of PWM actuation in the model (this algorithm is denoted as impulsive MPC). Finally, closed-loop simulations using MPC, based on the PWM algorithms as explained in Section~\ref{sec-mpcrz}. The impulses produced by the first and second methods are subsequently transformed to PWM inputs using the algorithm of Section~\ref{sec-pwmfilt}.

Compare first the algorithms  without disturbances. The trajectories (projected on the target orbital plane) are shown in Fig.~\ref{fig-traj}. The open-loop impulsive solution does not achieve rendezvous and drifts away, whereas the other solutions successfully reach the origin (the simulation is stopped when the chaser vehicle was 5 meters or less away from the target). The impulsive MPC is able to mostly compensate its imperfect thruster model. The PWM MPC  algorithm had a cost of $15.0~\mathrm{m/s}$ and the impulsive MPC had a cost of $15.8~\mathrm{m/s}$. Thus, while a basic MPC is able to rendezvous, the use of an imperfect model has some fuel costs. In addition, the impulsive MPC does not satisfy the line-of-sight constraints for a period of time.

Next, Fig.~\ref{fig-dist}  shows a simulation where the real orbit is different from the reference orbit used in the model (the real eccentricity is $e=0.83$, the real perigeee altitude is $h_p=525~\mathrm{km}$, and the real $\theta_0=60^\mathrm{o}$). Both MPC algorithms reach the origin (as in the previous scenario, the simulation is stopped when the chaser vehicle was 5 meters or less away from the target). The impulsive MPC again exits the line-of-sight region. The cost for the PWM MPC algorith was $15.3~\mathrm{m/s}$, whereas the impulsive MPC had a cost of $15.8~\mathrm{m/s}$.

Each iteration took less than half a second on a conventional computer, using MATLAB and the  \emph{Gurobi} optimization package (see~\cite{gurobi}). With a maximum number of iterations of 6, the computation time remained well below the interval sampling time.

\section{Concluding Remarks} \label{sect-rem}

This paper has presented a MPC algorithm that computes optimal PWM inputs for LTV systems. The algorithm is based on an initial approximation with either PAM or impulsive inputs, followed by iterative explicit linearization. As an application, the problem of rendezvous in elliptical orbits has been considered. In particular, the algorithm might be particularly useful for satellites with small specific thrust. The algorithm improves the fuel cost of an impulsive-only MPC (with the impulses posteriorly transformed to PWM inputs), and is able to satisfy safety constraints and handle disturbances such as imperfect knowledge of the target's orbit. This algorithm would help avoiding having to include a PWM approximation term  in the ``uncertainty budget'' and therefore save costs. However, inclusion of real-life constraints and more realistic simulations are needed to validate the method.

Possible future lines of research include studying the convergence of the planning algorithm, guaranteeing constraint satisfaction by including an estimate of linearization error in the model, or analyzing the stability guarantees of the MPC design.

\section*{Acknowledgments}

The authors acknowledge financial
support of the Spanish Ministry  of Science and Innovation under grant
DPI2008-05818.

\section*{References}
\bibliographystyle{elsarticle-harv.bst}
\bibliography{Rendezvous}
\newpage
\listoffigures

\end{document}